\documentclass{amsart}

\usepackage{amsmath}
\usepackage{amssymb}
\usepackage{amsthm}
\usepackage{xypic}
\usepackage{algorithmic}
\usepackage{algorithm}
\usepackage[pdftex]{graphicx}

\newtheorem{theorem}{Theorem}[section]
\newtheorem*{theorem*}{Theorem}
\newtheorem*{corollary*}{Corollary}
\newtheorem{lemma}[theorem]{Lemma}
\newtheorem{proposition}[theorem]{Proposition}

\newtheorem{corollary}[theorem]{Corollary}

\numberwithin{equation}{section}
\numberwithin{table}{section}

\def\Z{\mathbb{Z}}

\def\Q{\mathbb{Q}}
\def\C{\mathbb{C}}
\def\R{\mathbb{R}}
\def\OK{\mathcal{O}_K}

\def\qed{\raisebox{-.25ex}{\scalebox{.786}[1.272]{$\blacksquare$}}}

\makeatletter
\def\pmod#1{\allowbreak\mkern10mu({\operator@font mod}\,\,#1)}

\makeatother

\author{Kevin J. McGown}
\date{}
\title[Norm-Euclidean Galois fields and the GRH]{Norm-Euclidean Galois fields and the\\ Generalized Riemann Hypothesis}
\address{Department of Mathematics, University of California, San Diego, 9500 Gilman Drive, La Jolla, CA 92093}
\curraddr{Department of Mathematics, Oregon State University, 368 Kidder Hall, Corvallis, OR 97331}
\email{mcgownk@math.oregonstate.edu}

\begin{document}

\begin{abstract}
Assuming the Generalized Riemann Hypothesis (GRH), 
we show that
the norm-Euclidean Galois cubic fields
are exactly those with discriminant
$$
    \Delta=7^2,9^2,13^2,19^2,31^2,37^2,43^2,61^2,67^2,103^2,109^2,127^2,157^2
    \,.
$$  
A large part of the proof is in establishing the following
more general result:
Let $K$ be a Galois number field of odd prime degree $\ell$ and conductor $f$.  Assume the GRH for $\zeta_K(s)$.
If
$$
  38(\ell-1)^2(\log f)^6\log\log f<f
  \,,
$$  
then $K$ is not norm-Euclidean.
\end{abstract}

\maketitle


\vspace{-0.25in}

\section{Introduction}\label{S:intro}
Let $K$ be a number field with ring of integers $\OK$, and denote by $N=N_{K/\Q}$ the
absolute norm map.  For brevity, we will sometimes use the term field to mean a number field.
We call a number field $K$ norm-Euclidean if for every $\alpha,\beta\in\OK$, $\beta\neq 0$, there exists
$\gamma\in\OK$ such that $|N(\alpha-\gamma\beta)|<|N(\beta)|$.
In the quadratic setting, it is known that there are only finitely many norm-Euclidean fields and they have been identified;
namely, a number field of the form $K=\Q(\sqrt{d})$ with $d$ squarefree is norm-Euclidean if and only if
$$
  d=-1,-2,-3,-7,-11,2, 3, 5, 6, 7, 11, 13, 17, 19, 21, 29, 33, 37, 41, 57, 73
  \,.
$$

The main goal of this paper is to prove the following:
\begin{theorem}\label{T:cubic.GRH}
Assuming the GRH,
the norm-Euclidean Galois cubic fields
are exactly those with discriminant
$$
    \Delta=7^2,9^2,13^2,19^2,31^2,37^2,43^2,61^2,67^2,103^2,109^2,127^2,157^2
    \,.
$$  
\end{theorem}
For most of this paper, the reader may take the Generalized Riemann Hypothesis (GRH) to mean that
for every Dirichlet L-function $L(s,\chi)$, all the zeros of $L(s,\chi)$ in the critical strip $0<\Re(s)<1$
are on the critical line $\Re(s)=1/2$.\footnote{Actually, for Theorem~\ref{T:cubic.GRH}, it suffices to assume the Riemann Hypothesis (RH) and the GRH
for Dirichlet L-functions associated to cubic characters.}
The only exceptions will be when we explicitly state which function is being referred to ---
i.e., ``the GRH for $L(s,\chi)$'' or ``the GRH for $\zeta_K(s)$''.

Previously, Heilbronn (see~\cite{heilbronn:cubic}) showed that there are finitely many norm-Euclidean Galois cubic fields, but produced no upper
bound on the discriminant.  Godwin and Smith (see~\cite{godwin.smith:1993}) showed that the Galois cubic fields with $|\Delta|<10^8$ are exactly those listed in Theorem~\ref{T:cubic.GRH}
and were the first to give this list.  Lemmermeyer subsequently extended this result to show that
Godwin and Smith's list includes all fields with $|\Delta|<2.5\cdot 10^{11}$ (see~\cite{lemmermeyer:euclidean}).
Although it is the natural question, no one seems to have conjectured that this is the complete list;
however,
in light of Theorem~\ref{T:cubic.GRH}, this now seems like a very reasonable conjecture!

In a recent paper (see~\cite{mcgown:euclidean}) the author proved the
following unconditional result:
\begin{theorem}\label{T:cubic}
The
fields listed in Theorem~\ref{T:cubic.GRH} are norm-Euclidean, and
any remaining norm-Euclidean Galois cubic field
must have discriminant
$\Delta=f^2$ with $f\equiv 1\pmod{3}$
where $f$
is a prime in the interval $(10^{10},\,10^{70})$.
\end{theorem}
A large part of the proof of Theorems~\ref{T:cubic.GRH} and~\ref{T:cubic}
is in giving an upper bound on the discriminant for the class of fields in question.
Our technique works not only in the case of Galois cubic fields,
but for Galois fields of odd prime degree.
\begin{theorem}\label{T:DB}
  Let $\ell$ be an odd prime.  There exists a
  computable constant $C_\ell$ such that if $K$ is
  a Galois number field of odd prime degree $\ell$, conductor $f$,
  and discriminant $\Delta$, which is norm-Euclidean,
  then $f<C_\ell$ and $0<\Delta<C_\ell^{\ell-1}$.
\begin{table}[H]
\centering
\begin{tabular}{c c}
\begin{tabular}{| l | l |}
\hline
$\ell$ & $\;C_\ell$ \\
\hline
$3$ & $10^{70}$\\
$5$ & $10^{78}$\\
$7$ & $10^{82}$\\
$11$ & $10^{88}$\\
$13$ & $10^{89}$\\
$17$ & $10^{92}$\\
$19$ & $10^{94}$\\
$23$ & $10^{96}$\\
\hline
\end{tabular}
\qquad
\begin{tabular}{| l | l |}
\hline
$\ell$ & $\;C_\ell$ \\
\hline
$29$ & $10^{98}$\\
$31$ & $10^{99}$\\
$37$ & $10^{101}$\\
$41$ & $10^{102}$\\
$43$ & $10^{102}$\\
$47$ & $10^{103}$\\
$53$ & $10^{104}$\\
$59$ & $10^{105}$\\
\hline
\end{tabular}
\qquad
\begin{tabular}{| l | l |}
\hline
$\ell$ & $\;C_\ell$ \\
\hline
$61$ & $10^{106}$\\
$67$ & $10^{107}$\\
$71$ & $10^{107}$\\
$73$ & $10^{108}$\\
$79$ & $10^{108}$\\
$83$ & $10^{109}$\\
$89$ & $10^{109}$\\
$97$ & $10^{110}$\\
\hline
\end{tabular}
\end{tabular}
\vspace{0.5ex}
\caption{Values of $C_\ell$ for primes $\ell<100$\label{Table:Cl}}
\end{table}
\vspace{-2ex}
  \begin{table}[H]
\centering
\begin{tabular}{c c}
\begin{tabular}{| l | l |}
\hline
$\ell$ & \;$C_\ell$ \\
\hline
$3$ & $10^{11}$\\
$5$ & $10^{12}$\\
$7$ & $10^{13}$\\
$11$ & $10^{13}$\\
$13$ & $10^{14}$\\
$17$ & $10^{14}$\\
$19$ & $10^{14}$\\
$23$ & $10^{14}$\\
\hline
\end{tabular}
\qquad
\begin{tabular}{| l | l |}
\hline
$\ell$ & \;$C_\ell$ \\
\hline
$29$ & $10^{15}$\\
$31$ & $10^{15}$\\
$37$ & $10^{15}$\\
$41$ & $10^{15}$\\
$43$ & $10^{15}$\\
$47$ & $10^{15}$\\
$53$ & $10^{15}$\\
$59$ & $10^{15}$\\
\hline
\end{tabular}
\qquad
\begin{tabular}{| l | l |}
\hline
$\ell$ & \;$C_\ell$ \\
\hline
$61$ & $10^{15}$\\
$67$ & $10^{15}$\\
$71$ & $10^{16}$\\
$73$ & $10^{16}$\\
$79$ & $10^{16}$\\
$83$ & $10^{16}$\\
$89$ & $10^{16}$\\
$97$ & $10^{16}$\\
\hline
\end{tabular}
\end{tabular}
\vspace{0.5ex}
\caption{Values of $C_\ell$ for primes $\ell<100$, assuming the GRH} 
\label{Table:Cl2}
\end{table}
\end{theorem}
In~\cite{mcgown:euclidean}, the author proved Theorem~\ref{T:DB} and gave the constants in Table~\ref{Table:Cl}.
In this paper, we will show that under the GRH
these constants can be improved to those
given in Table~\ref{Table:Cl2}.  In fact, we will prove the following result
which,
after some easy computation,
completely justifies Table~\ref{Table:Cl2}.

\begin{theorem}\label{T:GRH}
 Let $K$ be a Galois number field of odd prime degree $\ell$ and conductor $f$.
 Assume the GRH for $\zeta_K(s)$, the Dedekind zeta function of $K$.\footnote{Note that, in this context, the GRH (as defined above) implies the GRH for $\zeta_K(s)$;
this follows from Lemma~\ref{L:factor}.
}
If
$$
  38(\ell-1)^2(\log f)^6\log\log f<f
  \,,
$$
then $K$ is not norm-Euclidean.
\end{theorem}


\section{Preliminaries}

As is customary, we will write $\zeta(s)$ to denote the Riemann zeta function,
$L(s,\chi)$ to denote the Dirichlet L-function associated to a Dirichlet character $\chi$,
and $\zeta_K(s)$ to denote the
Dedekind zeta function associated to a number field $K$.
The following is well-known and is an easy consequence of Theorem~8.6 of~\cite{narkiewicz:book}.

\begin{lemma}\label{L:factor}
Let $K$ be a Galois number field of odd prime degree $\ell$ and conductor $f$,
and  let $\chi$ be a primitive Dirichlet character modulo $f$ of order $\ell$.
Then
$$
  \zeta_K(s)=\zeta(s)\prod_{k=1}^{\ell-1}L(s,\chi^k)
  \,.
$$
\end{lemma}

We now quote three results from~\cite{mcgown:euclidean}
which will be crucial for our arguments.

\begin{lemma}\label{L:disc}
Let $K$ be a Galois number fields of odd prime degree $\ell$, conductor $f$, and discriminant $\Delta$.
Further, suppose that $K$ has class number one.
In this case, one has
$\Delta=f^{\ell-1}$.  Moreover:
\begin{enumerate}
\item
If $(f,\ell)=1$, then $f$ is a prime with $f\equiv1\pmod\ell$.
\item
If $(f,\ell)>1$, then $f=\ell^2$.
\end{enumerate}
\end{lemma}

\begin{theorem}\label{T:conditions}
  Let $K$ be a Galois number field of odd prime degree $\ell$ and conductor $f$ with $(f,\ell)=1$,
  and
  let $\chi$ be a primitive Dirichlet character modulo $f$ of order $\ell$.
  Denote by $q_1<q_ 2$ the two smallest rational primes that are inert in $K$.
  Suppose that there exists $r\in\Z^+$ with
  $$
    (r,q_1 q_2)=1,
    \quad
    \chi(r)=\chi(q_2)^{-1},
  $$
  such that any of the following conditions hold:
  \begin{enumerate}
    \item
    \qquad
    $r q_2 k\not\equiv f \pmod{q_1^2}$, \quad $k=1,\dots,q_1-1$,
    
    \quad\hspace{-1ex}
    $(q_1-1)(q_2 r-1)\leq f$
    
  \item
  \qquad
  $q_1\neq 2,3$,
  \quad
  $3q_1 q_2 r\log q_1< f $

  \item
  \qquad
  $q_1\neq 2,3,7$,
  \quad
  $2.1\,q_1 q_2 r\log q_1< f$
  
  \item
  \qquad
  $q_1=2,\,q_2\neq 3$,
  \quad
  $3q_2 r< f$
  
  \item
  \qquad
  $q_1=3,\,q_2\neq 5$,
  \quad
  $5q_2 r< f$
  \end{enumerate}
   Then $K$ is not norm-Euclidean.
\end{theorem}

\begin{proposition}\label{P:special}
Let $K$ be a Galois number field of odd prime degree $\ell$ and conductor $f$.
Denote by $q_1<q_ 2$ the two smallest rational primes that are inert in $K$.
Suppose either of the following conditions hold:
\begin{enumerate}
  \item
    \qquad
    $q_1=2,\, q_2=3$,
    
    \quad\hspace{-1ex}
    $72(\ell-1)f^{1/2}\log 4f+35\leq f$
  \item
    \qquad
      $q_1=3,\, q_2=5$,
      
    \quad\hspace{-1ex}
      $507(\ell-1)f^{1/2}\log 9f+448\leq f$
\end{enumerate}
Then $K$ is not norm-Euclidean.
\end{proposition}

\section{GRH Bounds for Non-Residues}\label{S:GRH.nonresidues}

In~\cite{bach:1990}, Bach proves an explicit version of a theorem due to Ankeny (see~\cite{ankeny}) regarding
the least element outside of a given non-trivial subgroup of $(\Z/m\Z)^\star$.  The main idea behind Bach's
proof appears in~\cite{montgomery:book}, but to obtain explicit results there are many details to work out;
Bach uses a slightly different kernel and introduces a parameter in order to achieve good numerical results.
Using the tables in~\cite{bach:1990}, we 
obtain the following special case which will be useful to us in the present context.
\begin{theorem}[Bach, 1990]\label{T:bach}
Assume the GRH.
Let $\chi$ be a non-principal Dirichlet character modulo $m\geq 10^8$,
and denote by $q_1$ the smallest prime such that
\mbox{$\chi(q_1)\neq 1$}.
Then
$$
  q_1<(1.17\,\log m-6.36)^2
  \,.
$$
\end{theorem}
We will follow Bach's approach to give bounds on $q_2$ and $r$.
Although the following results undoubtably hold in more generality,
we will not hesitate to specialize to our situation when it affords us certain technical conveniences.

\begin{theorem}\label{T:GRHq2}
Let $\chi$ be a non-principal Dirichlet character modulo $m\geq 10^{9}$ with $\chi(-1)=1$.
Assume the RH and the GRH for $L(s,\chi)$.
Denote by $q_1<q_2$ the two smallest primes such that $\chi(q_1),\chi(q_2)\neq 1$.
Then
$$
  q_2<2.5(\log m)^2
  \,.
$$
\end{theorem}

\begin{theorem}\label{T:GRHr}
Let $\ell$ and $f$ be odd primes with $f\geq 10^8$ and $f\equiv 1\pmod{\ell}$.
  Let $K$ be the Galois number field of degree $\ell$ and conductor $f$,
and let $\chi$ be a primitive Dirichlet character modulo $f$ of order $\ell$.
  Assume the GRH for $\zeta_K(s)$.
  If $q_1,q_2$ are rational primes and
  $\omega\neq 1$ is an $\ell$-th root of unity, then
  there exists $r\in\Z^+$ such that $(r,q_1 q_2)=1$, $\chi(r)=\omega$, and
  $$
    r<2.5(\ell-1)^2(\log f)^2
    \,.
  $$
\end{theorem}


The remainder of \S\ref{S:GRH.nonresidues} will be devoted to proving
Theorems~\ref{T:GRHq2} and~\ref{T:GRHr}.
Although this constitutes the bulk of the paper and is where the analytic techniques come into play,
the casual reader who is willing to accept these two results may skip the rest of this section and proceed to \S\ref{S:GRH.DB}.

In \S\ref{S:GRH.nonresidues.explicit} we give some explicit formulas
relating sums over prime powers to sums over zeros of L-functions,
and in \S\ref{S:GRH.nonresidues.zeros} we give some GRH estimates for
the sums over zeros.
Then in \S\ref{S:GRH.nonresidues.q2} and
\S\ref{S:GRH.nonresidues.r} we prove 
Theorems~\ref{T:GRHq2} and~\ref{T:GRHr}, respectively.


\subsection{An explicit formula}\label{S:GRH.nonresidues.explicit}

\begin{lemma}\label{L:explicit0}
Let $\chi$ be a Dirichlet character modulo $m$.
(Here we allow the possibility that $\chi$ is the principal character or even that $m=1$.)  For $x>1$ and $a\in(0,1)$, we have
$$
  -\frac{1}{2\pi i}
  \int_{2-i\infty}^{2+i\infty}
  \frac{x^s}{(s+a)^2}
  \frac{L'(s,\chi)}{L(s,\chi)}\,ds
  =
  \sum_{n<x}
  \chi(n)\Lambda(n)(n/x)^a\log(x/n)
  \,.
$$
\end{lemma}

\noindent\textbf{Proof.}
This is Lemma 4.2 of~\cite{bach:1990}.  We provide only a brief sketch here.
We plug the Dirichlet series
$$
  \frac{L'(s,\chi)}{L(s,\chi)}
  =
  -
  \sum_{n=1}^\infty
  \chi(n)\Lambda(n)n^{-s}
$$
into the left-hand side above and interchange the order of summation and integration.
Next, we use the fact that for $y>0$ one has
$$
  \frac{1}{2\pi i}
  \int_{2-i\infty}^{2+i\infty}
  \frac{y^s}{(s+a)^2}\;ds
  =
  \begin{cases}
  y^{-a}\log y & \text{ if $y>1$ }\\
  0 & \text{ otherwise }
  \end{cases}
  \,,
$$
and the result follows.
\qed

\begin{lemma}\label{L:explicit1}
  Let $\chi$ be a non-principal primitive Dirichlet character modulo $m$ with $\chi(-1)=1$.
  For $x>1$ and $a\in(0,1)$ we have
  \begin{eqnarray*}
  \sum_{n<x}
  \chi(n)\Lambda(n)(n/x)^a\log(x/n)
  &=&
  -\sum_{\rho\text{ of $L_\chi$}}
  \frac{x^\rho}{(\rho+a)^2}
  -
  \sum_{n=1}^\infty
  \frac{x^{-2n}}{(a-2n)^2}
    -
  \frac{1}{a^2}
  \\
  &&
  \quad
  -
  \frac{\log x}{x^a}
  \left(\frac{L_\chi'}{L_\chi}\right)(-a)
  -
  \frac{1}{x^a}
  \left(\frac{L_\chi'}{L_\chi}\right)'(-a)
  \,.
  \end{eqnarray*}
\end{lemma}

\noindent\textbf{Proof.}
Formally, this follows immediately by evaluating the integral in Lemma~\ref{L:explicit0} by residues.
For more details, see Lemma 4.4 of \cite{bach:1990}.
\qed

\begin{lemma}\label{L:explicit2}
  For $x>1$ and $a\in(0,1)$ we have
  \begin{eqnarray*}
  \sum_{n<x}
  \Lambda(n)(n/x)^a\log(x/n)
  &=&
   \frac{x}{(a+1)^2}
  -
  \sum_{\rho\text{ of $\zeta$}}
  \frac{x^\rho}{(\rho+a)^2}
   -
  \sum_{n=1}^\infty
  \frac{x^{-2n}}{(a-2n)^2}
  \\
  &&
  \quad
  -
  \frac{\log x}{x^a}
  \left(\frac{\zeta'}{\zeta}\right)(-a)
  -
  \frac{1}{x^a}
  \left(\frac{\zeta'}{\zeta}\right)'(-a)
  \,.
  \end{eqnarray*}
\end{lemma}

\noindent\textbf{Proof.}
This is similar to the previous result.
\qed

\vspace{1ex}

For our bounds on $q_2$ and $r$, we will need to exclude certain primes
from consideration; this will require the following estimate:
\begin{lemma}\label{L:sumbound}
 Let $u\in\Z^+$.  Then
 $$ 
  \sum_{\footnotesize{\begin{array}{cc}n<x\\[-.5ex](n,u)>1\end{array}}}
 \Lambda(n)
 (n/x)^a
 \log(x/n) 
 \leq
 \omega(u)(\log x)^2
 \,,
 $$
 where $\omega(u)$ denotes the number of distinct prime factors of $u$.
\end{lemma}

\noindent\textbf{Proof.}
If $u=1$, the result is trivial.
Suppose $u=p_1^{a_1}\dots p_t^{a_t}$.  Then
\begin{eqnarray*}
  \sum_{\footnotesize{\begin{array}{cc}n<x\\[-.5ex](n,u)>1\end{array}}}
  \Lambda(n)
   \;=\;
   \sum_{k=1}^t
   \sum_{a=1}^{\lfloor \log_{p_k}x\rfloor}
   \log p_k
   \;\leq\;
   \sum_{k=1}^t
   \log x
   \;=\;
   t\log x
   \,.
\end{eqnarray*}

The result easily follows.
\qed


\subsection{Sums over zeros}\label{S:GRH.nonresidues.zeros}
In order to prove our results, we will need to bound the sums over zeros appearing in
Lemmas~\ref{L:explicit1} and \ref{L:explicit2}.  Eventually we will take character combinations
of the formulas appearing in these lemmas as well, and so it will be useful to bound
the corresponding sum over all the zeros of the Dedekind zeta function of a number field $K$.

Let $K$ be a number field of discriminant $\Delta$ with $r_1$ real embeddings and $2r_2$ complex embeddings.  We define
$$
  \psi_K(s)
  :=
  \frac{r_1+r_2}{2}\,
  \psi\left(\frac{s}{2}\right)
  +
  \frac{r_2}{2}\,
  \psi\left(\frac{s+1}{2}\right)
  -
  \frac{1}{2}\,
  [K:\Q]\log \pi
  \,,\quad
  \psi(s):=\frac{\Gamma'(s)}{\Gamma(s)}
  \,,
$$
where $\Gamma(s)$ is the usual gamma function.
In particular,
$$
  \psi_\Q(s) = \frac{1}{2}\left(\psi\left(\frac{s}{2}\right)-\log\pi\right)
  \,.
$$
In order to expedite the proofs of this section, we quote some formulae,
all of which can be derived from (5.9) of~\cite{lagarias:1977}.
For all $s\in\C$, we have:
\begin{eqnarray}
 \label{E:zetaformula1}
  \frac{\zeta_K'(s)}{\zeta_K(s)}
  &=&
  B_K
  +
  \sum_{\rho\text{ of $\zeta_K$}}
  \left(
  \frac{1}{s-\rho}+\frac{1}{\rho}
  \right)
  -
  \frac{1}{2}
  \log|\Delta|
  -\frac{1}{s}-\frac{1}{s-1}-\psi_K(s)
  \\
   \label{E:zetaformula2}
  \frac{\zeta'(s)}{\zeta(s)}
  &=&
  B
  +
  \sum_{\rho\text{ of $\zeta$}}
  \left(
  \frac{1}{s-\rho}
  +
  \frac{1}{\rho}
  \right)
  -\frac{1}{s}
  -\frac{1}{s-1}
  -\psi_\Q(s)
\end{eqnarray}
If $\chi$ is a non-principal primitive Dirichlet character modulo $f$,
with $\chi(-1)=1$, then for all $s\in\C$ we have:
\begin{eqnarray}
   \label{E:zetaformula3}
    \frac{L'(s,\chi)}{L(s,\chi)}
  &=&
  B_\chi+
  \sum_{\rho\text{ of $L_\chi$}}
  \left(
  \frac{1}{s-\rho}
  +
  \frac{1}{\rho}
  \right)
  -\frac{1}{2}\log f
  -\psi_\Q(s)
\end{eqnarray}
Each sum above is over the non-trivial zeros $\rho$ of the corresponding functions,
and is absolutely and uniformly convergent on compact subsets of $\C$.
Henceforth we adopt the notation that $\rho$ will always denote a non-trivial zero
with $0<\Re(\rho)<1$.

Each of (\ref{E:zetaformula1}), (\ref{E:zetaformula2}), (\ref{E:zetaformula3}) involves a constant $B$ which can be difficult to estimate.
Fortunately, in all three cases this constant can be eliminated from the equation
as follows.
Provided the sum is taken in symmetric order\footnote
{
Taking the sum in symmetric order means:  $\;\displaystyle\sum_\rho \;=\; \lim_{T\to\infty}\sum_{\substack{\rho=\sigma+it\\|t|<T}}$
}
,
one has
\begin{eqnarray}
  &&
  B+\sum_{\rho\text{ of $\zeta$ }}\frac{1}{\rho}=0
  \,,
\end{eqnarray}
and similarly for $B_K$ and $B_\chi$.
See~\cite{davenport:mnt} for a simple argument which gives this result for the constant $B$.
The corresponding result for $B_K$ follows by a similar argument and
was first exploited by Stark to give lower bounds for discriminants (see~\cite{stark:1974, stark:1975}).

The analogous result for $B_\chi$
is not obvious; in fact, it wasn't known until the
introduction of the Weil formulas (see~\cite{weil:ef.1, weil:ef.2}).
Plugging $s=1$ into (\ref{E:zetaformula3}) and comparing against (2.3.1) of~\cite{murty:2009} 
gives a proof of this result.
See \cite{poitou:1, poitou:2, odlyzko:1990} for results regarding the use of explicit formulae
to obtain discriminant bounds.

We begin with a lemma which goes back to Landau (see~\cite{landau:1}).
\begin{lemma}\label{L:stark}
Let $\chi$ be a primitive Dirichlet character modulo $f$ with $\chi(-1)=1$.
Then for $\sigma\in\R$, we have
$$
\sum_{\rho\text{ of } L_\chi}
\left(
\frac{1}{\sigma-\rho}
+
\frac{1}{\sigma-\overline{\rho}}
\right)
=
\log f
+
2
\Re
\frac{L'(\sigma,\chi)}{L(\sigma,\chi)}
+
2\psi_\Q(\sigma)
\,.
$$
\end{lemma}

\noindent\textbf{Proof.}
We substitute $s=\sigma$ into
(\ref{E:zetaformula3})
and add the result to it's conjugate.
The result now follows upon invoking the fact that
$$
  \Re\left(
  B_\chi
  +
  \sum_{\rho\text{ of } L_\chi}
  \frac{1}{\rho}
  \right)
  =0
  \,.
  \;\;
  \text{\qed}
$$ 

\begin{lemma}\label{L:justlike}
Let $\chi$ be a non-principal primitive Dirichlet character modulo $f$ with $\chi(-1)=1$.
Assume the RH and the GRH for $L(s,\chi)$.
For $a\in(0,1)$ we have
$$
    \sum_{\rho\text{ of $\zeta,L_\chi$ }}
  \frac{1}{|\rho+a|^2}
  \leq
  \frac{1}{2a+1}
  \left(
  \log f
  +
  2
  \left(\frac{1}{a+1}+\frac{1}{a}\right)
  +
  4\psi_\Q(a+1)
  \right)
  \,.
$$
\end{lemma}

\noindent\textbf{Proof.}
We consider the following two formulae:
\begin{eqnarray*}
\sum_{\rho\text{ of } \zeta}
\left(
\frac{1}{\sigma-\rho}
+
\frac{1}{\sigma-\overline{\rho}}
\right)
&=&
2\frac{\zeta'(\sigma)}{\zeta(\sigma)}
+
2\left(\frac{1}{\sigma}+\frac{1}{\sigma-1}\right)
+
2\psi_\Q(\sigma)
\\[1ex]
\sum_{\rho\text{ of } L_\chi}
\left(
\frac{1}{\sigma-\rho}
+
\frac{1}{\sigma-\overline{\rho}}
\right)
&=&
\log f
+
2
\Re
\frac{L'(\sigma,\chi)}{L(\sigma,\chi)}
+
2\psi_\Q(\sigma)
\end{eqnarray*}
The second formula above
is Lemma~\ref{L:stark} and the first can be proved in exactly the same manner.
Setting $\sigma=a+1$ and supposing that
$\Re(\rho)=1/2$, we find:
\begin{equation}\label{E:GRH.alg}
  \frac{1}{|\rho+a|^2}
  =
  \frac{1}{2a+1}
  \left(
  \frac{1}{\sigma-\rho}+\frac{1}{\sigma-\overline{\rho}}
  \right)
\end{equation}
To complete the proof, we combine everything above and note that
$$
  \frac{\zeta'(\sigma)}{\zeta(\sigma)}
  +
  \Re
  \frac{L'(\sigma,\chi)}{L(\sigma,\chi)}
  <
  0
  \,,
$$
by considering the Dirichlet series for $(\zeta'/\zeta+L_\chi'/L_\chi)(s)$.
\qed

\vspace{1ex}

We give a special case of the previous lemma:
\begin{lemma}\label{L:q2zerosbound}
Let $\chi$ be a non-principal primitive Dirichlet character modulo $f$ with $\chi(-1)=1$.
Assume the RH and the GRH for $L(s,\chi)$.  We have
$$
  \sum_{\rho\text{ of $\zeta,L_\chi$ }}
  \frac{1}{\left|\rho+\frac{1}{2}\right|^2}
  \leq
  \frac{1}{2}\log f
  +0.437
  \,.
$$

\end{lemma}

\noindent\textbf{Proof.}
Use the fact
\begin{equation}\label{E:somecomp}
  \psi_\Q(3/2)\approx
  -1.1153
\end{equation}
and apply the previous lemma with $a=1/2$.
\qed

\vspace{1ex}

Having completed the desired estimates over the zeros of $\zeta(s)$ and $L(s,\chi)$, we turn
turn to $\zeta_K(s)$.

\begin{lemma}\label{L:stark2}
Let $K$ be a number field with discriminant $\Delta$.
Then we have
$$
\sum_{\rho\text{ of $\zeta_K$}}
\left(
\frac{1}{\sigma-\rho}+\frac{1}{\sigma-\overline{\rho}}
\right)
=
\log|\Delta|
+
2
\left(\frac{1}{\sigma}+\frac{1}{\sigma-1}\right)
+
2\psi_K(\sigma)
+
2
\frac{\zeta_K'(\sigma)}{\zeta_K(\sigma)}
\,.
$$
\end{lemma}

\noindent\textbf{Proof.}
This is exactly analogous to Lemma~\ref{L:stark}.
\qed

\begin{lemma}
Let $K$ be a number field with discriminant $\Delta$.
Suppose the GRH holds for $\zeta_K(s)$.
For $a\in(0,1)$ we have
\begin{eqnarray*}
  \sum_{\rho\text{ of $\zeta_K$}}
  \frac{1}{|\rho+a|^2}
  &<&
  \frac{1}{2a+1}
  \left[
  \log|\Delta|
  +
  2
  \left(\frac{1}{a+1}+\frac{1}{a}\right)
  +
  2
  \psi_K(a+1)
  \right]
  \,.
\end{eqnarray*}
\end{lemma}

\noindent\textbf{Proof.}
Let $\sigma=a+1$.
%
Applying Lemma~\ref{L:stark2} and using (\ref{E:GRH.alg}) gives
$$
  \sum_{\rho\text{ of $\zeta_K$}}
  \frac{1}{|\rho+a|^2}
  =
  \frac{1}{2a+1}
  \left[
  \log|\Delta|
  +
  2
  \left(\frac{1}{a+1}+\frac{1}{a}\right)
  +
  2
  \psi_K(a+1)
    +
  2\,\frac{\zeta_K'(a+1)}{\zeta_K(a+1)}
  \right]
  \,.
$$
The result follows upon observing that $\zeta_K'(\sigma)/\zeta_K(\sigma)<0$.
\qed

\vspace{1ex}

We give a special case of the previous lemma:
\begin{lemma}
Let $K$ be a totally real number field with discriminant $\Delta$.
Suppose the GRH holds for $\zeta_K(s)$.
We have
$$
  \sum_{\rho\text{ of $\zeta_K$}}
  \frac{1}{\left|\rho+\frac{1}{2}\right|^2}
<
\frac{1}{2}
\left(
\log|\Delta|
+
\frac{16}{3}
+
2\,\psi_\Q(3/2)\,
[K:\Q]
\right)
\,.
$$
\end{lemma}

\noindent\textbf{Proof.}
Since $r_1=[K:\Q]$, $r_2=0$, we have
$
\psi_K(s)=[K:\Q]\psi_\Q(s)
\,.
$
The result now follows from the previous lemma
upon setting $a=1/2$.
\qed

\vspace{1ex}

Now we specialize even further to our situation:
\begin{lemma}\label{L:rzerosbound}
Let $K$ be a totally real number field of degree $\ell$ and discriminant $\Delta=f^{\ell-1}$.
Suppose the GRH holds for $\zeta_K(s)$.
We have
$$
  \sum_{\rho\text{ of $\zeta_K$}}
  \frac{1}{\left|\rho+\frac{1}{2}\right|^2}
<
  \frac{1}{2}
  \left[
  (\ell-1)\log f-
  2.23\,\ell
  +5.34
  \right]
  \,.
$$
\end{lemma}

\noindent\textbf{Proof.}
We apply the previous lemma, using the approximation
given in (\ref{E:somecomp}).
\qed
%


\subsection{An upper estimate on $q_2$}\label{S:GRH.nonresidues.q2}

We establish a series of results, building up to the proof of Theorem~\ref{T:GRHq2}.

\begin{lemma}\label{L:q2lemma1}
Let $\chi$ be a non-principal Dirichlet character modulo $m$ with \mbox{$\chi(-1)=1$}.
For $a\in(0,1)$ and $x>0$ we have
\begin{eqnarray*}
&&
\frac{x}{(a+1)^2}
+
\frac{1}{a^2}
 \;\;=\;\;
  \sum_{\rho\text{ of $\zeta$ }}
   \frac{x^\rho}{(\rho+a)^2}
  -
   \sum_{\rho\text{ of $L_\chi$ }}
 \frac{x^\rho}{(\rho+a)^2}
 \\ 
 &&\qquad\qquad\qquad\qquad
  \;+\;
    \sum_{\footnotesize{\begin{array}{cc}n<x\\[-1ex]\chi(n)\neq 1\end{array}}}
    (1-\chi(n))
 \Lambda(n)
 (n/x)^a
 \log(x/n) 
 \\ 
 &&\qquad\qquad\qquad\qquad
  \;+\;
 \frac{\log x}{x^a}
  \left[
  \left(\frac{\zeta'}{\zeta}\right)(-a)
  -
 \left(\frac{L_{\chi}'}{L_{\chi}}\right)(-a)
 \right]
 \\
 &&\qquad\qquad\qquad\qquad
 \;+\;
 \frac{1}{x^a}
 \left[
  \left(\frac{\zeta'}{\zeta}\right)'(-a)
  -
 \left(
 \frac{L_{\chi}'}{L_{\chi}}
 \right)'
 (-a)
 \right]
 \,.
\end{eqnarray*}
\end{lemma}

\noindent\textbf{Proof.}
Subtract Lemma~\ref{L:explicit1} from Lemma~\ref{L:explicit2}.
\qed

\begin{lemma}\label{L:q2lemma2}
Let $\chi$ be a non-principal primitive Dirichlet character modulo $f$ with \mbox{$\chi(-1)=1$}.
For $a\in(0,1)$ we have
\begin{eqnarray*}
   &&
    \left|
  \left(\frac{\zeta'}{\zeta}\right)(-a)
  -
  \left(\frac{L_{\chi}'}{L_{\chi}}\right)(-a)
  \right|
  \\[1ex]
  &&
  \qquad\qquad
  \leq
  (a+2)
   \sum_{\rho\text{ of $\zeta,L_\chi$ }}
  \frac{1}{|(\rho+a)(2-\rho)|}
  +
  2\,
  \left|\frac{\zeta'(2)}{\zeta(2)}\right|
  +
  \frac{1}{a}
  +
  \frac{1}{a+1}
  +
  \frac{3}{2}
  \,.
\end{eqnarray*} 
\end{lemma}

\noindent\textbf{Proof.}
We begin with the following formulas which hold for all $s\in\C$, provided the sums are taken in symmetric order:
\begin{eqnarray}
  \label{E:Lformula1}
  \left(\frac{\zeta'}{\zeta}\right)(s)
  &=&
   \sum_{\rho\text{ of $\zeta$ }}
  \frac{1}{s-\rho}
  -\frac{1}{s}
  -\frac{1}{s-1}
  -\psi_\Q(s)
  \\
  \label{E:Lformula2}
  \left(\frac{L_\chi'}{L_\chi}\right)(s)
  &=&
   \sum_{\rho\text{ of $L_\chi$ }}
  \frac{1}{s-\rho}
  -\frac{1}{2}\log f
  -\psi_\Q(s)
\end{eqnarray}
Formulas (\ref{E:Lformula1}) and (\ref{E:Lformula2}) are obtained from
(\ref{E:zetaformula2}) and (\ref{E:zetaformula3}) respectively by applying the facts
$\sum_{\rho\text{ of $\zeta$ }}\rho^{-1}+B=0$ and $\sum_{\rho\text{ of $L_\chi$ }}\rho^{-1}+B_\chi=0$.
Plugging $s=2$ into (\ref{E:Lformula1}) and subtracting it from itself, and similarly for (\ref{E:Lformula2}), yields:
\begin{eqnarray*}
&&
\hspace{-2ex}
  \left(\frac{\zeta'}{\zeta}\right)(s)
  =
  \left(\frac{\zeta'}{\zeta}\right)(2)
  +
  \sum_\rho
  \left(\frac{1}{s-\rho}-\frac{1}{2-\rho}\right)
  +\frac{3}{2}
  -\frac{1}{s}
  -\frac{1}{s-1}
  +
  \psi_\Q(2)-\psi_\Q(s)
  \\
  &&
  \hspace{-2.75ex}
  \left(\frac{L_\chi'}{L_\chi}\right)(s)
=
  \left(\frac{L_\chi'}{L_\chi}\right)(2)
  +
  \sum_\rho
    \left(\frac{1}{s-\rho}-\frac{1}{2-\rho}\right)
  +
    \psi_\Q(2)-\psi_\Q(s)
\end{eqnarray*}
Using the above, together with the fact
$$
 \frac{1}{-a-\rho}-\frac{1}{2-\rho}
 =
 -
 \frac{a+2}{(\rho+a)(2-\rho)}
 \,,
$$
we can write
\begin{eqnarray*}
  &&
  \left(\frac{\zeta'}{\zeta}\right)(-a)
  -
  \left(\frac{L_{\chi}'}{L_{\chi}}\right)(-a)
  \\[.5ex]
  &&
  \qquad
  =
  (a+2)
  \left(
  \sum_{\rho\text{ of $L_\chi$}}
   \frac{1}{(\rho+a)(2-\rho)}
  -
  \sum_{\rho\text{ of $\zeta$}}
   \frac{1}{(\rho+a)(2-\rho)}
   \right)
  \\[.5ex]
  &&
  \qquad
  \qquad
  +
  \left(\frac{\zeta'}{\zeta}\right)(2)
  -
  \left(\frac{L_\chi'}{L_\chi}\right)(2)
    +
  \frac{3}{2}+\frac{1}{a}+\frac{1}{a+1}
  \,.
\end{eqnarray*}
The result follows upon taking absolute values and using the fact that
$$
  \left|
  \left(
  \frac{L'_\chi}{L_\chi}  
  \right)
  (2)
  \right|
  \leq
    \left|
  \left(
  \frac{\zeta'}{\zeta}  
  \right)
  (2)
  \right|
  \,.
  \;\;
  \text{\qed}
$$

\begin{lemma}\label{L:q2lemma3}
Suppose $a\in(0,1)$ and $\Re(\rho)=1/2$.
Then
$$
  \frac{1}{|(\rho+a)(2-\rho)|}
  \leq
  \frac{1}{|\rho+a|^2}
  \,.
$$
\end{lemma}

\noindent\textbf{Proof.}
Use $|2-\rho|\geq|\rho+a|$.
\qed

\begin{lemma}\label{L:q2lemma4}
Let $\chi$ be a non-principal primitive Dirichlet character modulo $f$ with \mbox{$\chi(-1)=1$}.
For $a\in(0,1)$ we have
$$
 \left|
  \left(\frac{\zeta'}{\zeta}\right)'(-a)
  -
 \left(
 \frac{L_{\chi}'}{L_{\chi}}
 \right)'
 (-a)
 \right|
 <
 \sum_{\rho\text{ of $\zeta,L_\chi$ }}
 \frac{1}{|\rho+a|^2}
 +
 \frac{1}{a^2}
 +
 \frac{1}{(a+1)^2}
 \,.
$$
\end{lemma}

\noindent\textbf{Proof.}
We start by differentiating (\ref{E:Lformula1}) and (\ref{E:Lformula2}); this gives
\begin{eqnarray}
  \label{E:L1formula1}
  \left(\frac{\zeta'}{\zeta}\right)'(s)
  &=&
  -\sum_{\rho\text{ of $\zeta$}}
  \frac{1}{(s-\rho)^2}
  +
  \frac{1}{s^2}
  +
  \frac{1}{(s-1)^2}
  -
  \psi_\Q'(s)
  \\
  \label{E:L1formula2}
  \left(\frac{L_\chi'}{L_\chi}\right)'(s)
  &=&
  -\sum_{\rho\text{ of $L_\chi$}}
  \frac{1}{(s-\rho)^2}
  -
  \psi_\Q'(s)
  \,,
\end{eqnarray}
which allows us to write
\begin{eqnarray*}
  &&
  \hspace{-1.5ex}
  \left(\frac{\zeta'}{\zeta}\right)'(-a)
  -
 \left(
 \frac{L_{\chi}'}{L_{\chi}}
 \right)'
 (-a)
 =
 \sum_{\rho\text{ of $L_\chi$ }}
 \frac{1}{(\rho+a)^2}
 -
 \sum_{\rho\text{ of $\zeta$ }}
 \frac{1}{(\rho+a)^2}
 +
 \frac{1}{a^2}
 +
 \frac{1}{(a+1)^2}
 \,.
\end{eqnarray*}
The result follows.
\qed

\begin{proposition}\label{P:q2}
Let $\chi$ be a non-principal primitive Dirichlet character modulo $f$ with $\chi(-1)=1$.
Assume the RH and the GRH for $L(s,\chi)$.
We define
$$
  \sum_\rho:=\sum_{\rho\text{ of $\zeta,L_\chi$}}
  \frac{1}{\left|\rho+\frac{1}{2}\right|^2}
  \,.
$$
For $x>0$ we have
\begin{eqnarray*}
\frac{x}{9/4}+4
&\leq&
\sqrt{x}\sum_\rho
\;+\;
2
\sum_{\footnotesize{\begin{array}{cc}n<x\\[-1ex]\chi(n)\neq 1\end{array}}}
\Lambda(n)(n/x)^{1/2}\log(x/n)
\\[1ex]
&&
\qquad
+
\frac{\log x}{\sqrt{x}}
\left(
\frac{5}{2}
\sum_\rho
+
\,
2
\,
\left|\frac{\zeta'(2)}{\zeta(2)}\right|
+
\frac{25}{6}
\right)
+
\frac{1}{\sqrt{x}}
\left(
\sum_\rho
+\,
\frac{40}{9}
\right)
\,.
\end{eqnarray*}

\end{proposition}

\noindent\textbf{Proof.}
Set $a=1/2$.  Combine Lemmas
\ref{L:q2lemma1}, \ref{L:q2lemma2}, \ref{L:q2lemma3}, and \ref{L:q2lemma4}.
\qed

\vspace{1ex}
\noindent\textbf{Proof of Theorem~\ref{T:GRHq2}.}
The result for a general character follows from the corresponding result for
primitive characters and hence we may assume $\chi$ is a primitive
character modulo $f$.

Define $x:=2.5(\log f)^2$.
Since $f\geq 10^{9}$, we have $x>1073$.
By way of contradiction,
suppose that $\chi(n)=1$ for all $n<x$ with $(n, q_1)=1$.
Under this assumption, we apply
Lemma~\ref{L:sumbound} with $u=q_1$, which gives
$$
  \sum_{\footnotesize{\begin{array}{cc}n<x\\[-1ex]\chi(n)\neq 1\end{array}}}
  \Lambda(n)(n/x)^{1/2}\log(x/n)
  \leq (\log x)^2
\,.
$$
Combining the above with Proposition~\ref{P:q2} and dividing by $\sqrt{x}$ yields
$$
  \frac{\sqrt{x}}{9/4}
  +
  \frac{4}{\sqrt{x}}
  \leq
  \sum_\rho
  +
  \frac{2(\log x)^2}{\sqrt{x}}
  +
  \frac{\log x}{x}
  \left[
  \frac{5}{2}\sum_\rho+
  2
  \left|\frac{\zeta'(2)}{\zeta(2)}\right|+\frac{25}{6}
  \right]
  +
  \frac{1}{x}
  \left[
  \sum_\rho
  +
  \frac{40}{9}
  \right]
  \,.
$$

By Lemma~\ref{L:q2zerosbound},
we have
$$
  \sum_\rho\leq
  \frac{1}{2}
  \log f
  +0.437
  \,,
$$
and in particular,
$$
  \frac{1}{\sqrt{x}}\sum_\rho\leq
  \frac{1}{3}
  \,.
$$
We see
\begin{eqnarray*}
  \frac{\log x}{x}
  \left[
  \frac{25}{6}
  \right]
  +
  \frac{1}{x}
  \left[
  \sum_\rho
  +
  \frac{40}{9}
  \right]
  \;\leq\;
  \frac{1}{\sqrt{x}}
  \left(
        \frac{\log x}{\sqrt{x}}
  \cdot
  \frac{25}{6}
  +
  \frac{1}{3}
  +
  \frac{1}{\sqrt{x}}
  \frac{40}{9}
  \right)
  \;<\;
  \frac{4}{\sqrt{x}}
  \,,
\end{eqnarray*}
and therefore
$$
  \frac{\sqrt{x}}{9/4}
  \leq
  \sum_\rho
  +
  \frac{2(\log x)^2}{\sqrt{x}}
  +
  \frac{\log x}{x}
  \left[
  \frac{5}{2}\sum_\rho+2
  \left|\frac{\zeta'(2)}{\zeta(2)}\right|
  \right]
  \,.
$$
We have
$$
  \frac{1}{\sqrt{x}}\left[
  \frac{5}{2}\sum_\rho
  +
  2\left|\frac{\zeta'(2)}{\zeta(2)}\right|
  \right]
  <
  0.869
$$
and
$$
\frac{\log x}{\sqrt{x}}
\leq
0.214
$$
which leads to
$$
  \frac{\sqrt{x}}{9/4}
  \leq
    \frac{1}{2}
    \log f
    +
    0.437
  +
  \frac{2(\log x)^2}{\sqrt{x}}
  +
  0.186
  \,.
$$
Now we observe
\begin{eqnarray*}
  \frac{2(\log x)^2}{\sqrt{x}}
  &\leq&
  2.98
  \,.
\end{eqnarray*}

All together, we have
\begin{eqnarray*}
  \frac{\sqrt{x}}{9/4}
  &\leq&
      \frac{1}{2}
      \log f
      +3.61
\end{eqnarray*}
This leads to:
\begin{eqnarray*}
  \sqrt{x}
  &\leq&
  \frac{9}{8}(\log f)
  +
  8.13
  \\
  &\leq&
  1.52\,\log f
\end{eqnarray*}
Squaring both sides yields
$$
  x\leq 2.32\,(\log f)^2
  \,,
$$  
a contradiction.
\qed


\subsection{An upper estimate on $r$}\label{S:GRH.nonresidues.r}

We establish a series of results, building up to the proof
of Theorem~\ref{T:GRHr}
\begin{lemma}\label{L:rlemma1}
Let $\chi$ be a non-principal Dirichlet character modulo a prime $p$ of order $\ell$ with $\chi(-1)=1$.
Fix any $\ell$-th root of unity $\omega\neq 1$.
For $a\in(0,1)$ and $x\in(1,p)$ we have
\begin{eqnarray*}
  &&
\frac{x}{(a+1)^2}
+
\frac{1}{a^2}
 \;\;=\;\;
 \sum_{k=1}^\ell
 \omega^{-k}
 \sum_{\rho\text{ of $L_{\chi^k}$}}
 \frac{x^\rho}{(\rho+a)^2}
 \\ 
 &&\qquad\qquad\qquad\qquad
  \;+\;
 \ell
  \sum_{\footnotesize{\begin{array}{cc}n<x\\[-1ex]\chi(n)=\omega\end{array}}}
 \Lambda(n)
 (n/x)^a
 \log(x/n) 
\\
 &&\qquad\qquad\qquad\qquad
 \;+\;
 \frac{\log x}{x^a}
  \left[
  \left(\frac{\zeta'}{\zeta}\right)(-a)
  +
   \sum_{k=1}^{\ell-1}
 \omega^{-k}
 \left(\frac{L_{\chi^k}'}{L_{\chi^k}}\right)(-a)
 \right]
 \\
 &&\qquad\qquad\qquad\qquad
 \;+\;
 \frac{1}{x^a}
 \left[
  \left(\frac{\zeta'}{\zeta}\right)'(-a)
  +
   \sum_{k=1}^{\ell-1}
 \omega^{-k}
 \left(
 \frac{L_{\chi^k}'}{L_{\chi^k}}
 \right)'
 (-a)
 \right]
\,.
\end{eqnarray*}
\end{lemma}

\noindent\textbf{Proof.}
First we note that $\chi^k$ for $k=1,\dots,\ell-1$ are all non-principal primitive characters
as $\chi$ is a character modulo a prime $p$ of order $\ell$;
moreover, $\chi^\ell(n)=1$ for all $n<x$ as $x<p$.
Multiplying the identity
$$
  \sum_{k=1}^\ell
  \omega^{-k}\chi^k(n)
  =
  \begin{cases}
  \ell & \chi(n)=\omega\\
  0 & \text{ otherwise }
  \end{cases}
$$
by
$$
  g(x,n):=\Lambda(n)(n/x)^a\log(x/n)
$$
and summing over all $n<x$ yields
$$
  \sum_{n<x}
  g(x,n)
  \sum_{k=1}^\ell
  \omega^{-k}
  \chi^k(n)
  =
  \ell
  \sum_{\footnotesize{\begin{array}{cc}n<x\\[-1ex]\chi(n)=\omega\end{array}}}
  g(x,n)
  \,.
$$
Interchanging the order of summation gives
$$
  \sum_{k=1}^\ell
  \omega^{-k}
  \sum_{n<x}
  g(x,n)\chi^k(n)
  =
  \ell
  \sum_{\footnotesize{\begin{array}{cc}n<x\\[-1ex]\chi(n)=\omega\end{array}}}
  g(x,n)
  \,.
$$
Now we apply Lemma~\ref{L:explicit1} and Lemma~\ref{L:explicit2} and use the facts:
\begin{equation}\label{E:thefacts}
  \sum_{k=1}^\ell
  \omega^{-k}
  =0
  \,,
  \quad
  \sum_{k=1}^{\ell-1}
  \omega^{-k}
  =-1
  \,.
\end{equation}
The result follows.
\qed

\begin{lemma}\label{L:rlemma2}
Let $\chi$ be a non-principal Dirichlet character modulo a prime $p$ of order $\ell$ with $\chi(-1)=1$.
Fix any $\ell$-th root of unity $\omega\neq 1$.
For $a\in(0,1)$ we have
\begin{eqnarray*}
   &&
    \left|
  \left(\frac{\zeta'}{\zeta}\right)(-a)
  +
   \sum_{k=1}^{\ell-1}
  \omega^{-k}
  \left(\frac{L_{\chi^k}'}{L_{\chi^k}}\right)(-a)
  \right|
  \\[1ex]
  &&
  \qquad\qquad
  \leq
  (a+2)
  \sum_\rho
  \frac{1}{|(\rho+a)(2-\rho)|}
  +
  \ell\,
  \left|\frac{\zeta'(2)}{\zeta(2)}\right|
  +
  \frac{1}{a}
  +
  \frac{1}{a+1}
  +
  \frac{3}{2}
  \,,
\end{eqnarray*} 
where the sum is taken over all non-trivial zeros $\rho$ of $L(s,\chi^k)$ for $k=1,\dots,\ell$. 
\end{lemma}

\noindent\textbf{Proof.}
Using
(\ref{E:Lformula1}), (\ref{E:Lformula2}) and (\ref{E:thefacts}),
we can write:
\begin{eqnarray*}
  &&
  \left(\frac{\zeta'}{\zeta}\right)(-a)
  +
   \sum_{k=1}^{\ell-1}
  \omega^{-k}
  \left(\frac{L_{\chi^k}'}{L_{\chi^k}}\right)(-a)
  \\
  &&
  \qquad
  =\;\;
    -(a+2)
     \sum_{k=1}^\ell
     \omega^{-k}
     \sum_{\rho\text{ of $L_{\chi^k}$ }}
     \frac{1}{(\rho+a)(2-\rho)}
   \\
   &&
   \qquad\qquad
     +
     \left(\frac{\zeta'}{\zeta}\right)(2)
    +
    \sum_{k=1}^{\ell-1}
      \left(\frac{L_{\chi^k}'}{L_{\chi^k}}\right)(2)
     +
     \frac{3}{2} 
     +
     \frac{1}{a}
     +
     \frac{1}{a+1}
\end{eqnarray*}
The result follows in a similar manner as Lemma~\ref{L:q2lemma2}.
\qed

%

\begin{lemma}\label{L:rlemma4}
Let $\chi$ be a non-principal Dirichlet character modulo a prime $p$ of order $\ell$ with $\chi(-1)=1$.
Fix any $\ell$-th root of unity $\omega\neq 1$.
For $a\in(0,1)$ we have
$$
 \left|
  \left(\frac{\zeta'}{\zeta}\right)'(-a)
  +
   \sum_{k=1}^{\ell-1}
 \omega^{-k}
 \left(
 \frac{L_{\chi^k}'}{L_{\chi^k}}
 \right)'
 (-a)
 \right|
 \leq
 \sum_{\rho}
 \frac{1}{|\rho+a|^2}
 +
 \frac{1}{a^2}
 +
 \frac{1}{(a+1)^2}
 \,,
$$
where the sum is taken over all non-trivial zeros $\rho$ of $L(s,\chi^k)$ for $k=1,\dots,\ell$. 
\end{lemma}

\noindent\textbf{Proof.}
Using (\ref{E:L1formula1}), (\ref{E:L1formula2}), and (\ref{E:thefacts}) we can write
\begin{eqnarray*}
  &&
  \left(\frac{\zeta'}{\zeta}\right)'(-a)
  +
  \sum_{k=1}^{\ell-1}
  \omega^{-k}
    \left(\frac{L_{\chi^k}'}{L_{\chi^k}}\right)'(-a)
    \\
    &&\qquad
    =\;
        \frac{1}{a^2}
  +
  \frac{1}{(a+1)^2}
    -\sum_{k=1}^\ell
    \omega^{-k}
    \sum_{\rho\text{ of $L_{\chi^k}$}}
      \frac{1}{(\rho+a)^2}
  \,.
  \;\;
  \text{\qed}
\end{eqnarray*}

%

\begin{proposition}\label{P:explicit}
Let $\ell$ and $f$ be odd primes with $f\equiv 1\pmod{\ell}$.
 Let $K$ be the Galois number field of degree $\ell$ and conductor $f$,
and let $\chi$ be a primitive Dirichlet character modulo $f$ of order $\ell$.
Fix any $\ell$-th root of unity $\omega\neq 1$.
Suppose that the GRH holds for $\zeta_K(s)$.
We define
$$
  \sum_\rho:=\sum_\rho
  \frac{1}{\left|\rho+\frac{1}{2}\right|^2}
  \,,
$$
where the sum is taken over all non-trivial zeros of $\zeta_K(s)$.
For $x\in(1,f)$ we have
\begin{eqnarray*}
\frac{x}{9/4}+4
&\leq&
\sqrt{x}\sum_\rho
\;+\;
\ell
\sum_{\footnotesize{\begin{array}{cc}n<x\\[-1ex]\chi(n)=\omega\end{array}}}
\Lambda(n)(n/x)^{1/2}\log(x/n)
\\[1ex]
&&
\qquad
+
\frac{\log x}{\sqrt{x}}
\left(
\frac{5}{2}
\sum_\rho
+
\,
\ell
\,
\left|\frac{\zeta'(2)}{\zeta(2)}\right|
+
\frac{25}{6}
\right)
+
\frac{1}{\sqrt{x}}
\left(
\sum_\rho
+\,
\frac{40}{9}
\right)
\,.
\end{eqnarray*}

\end{proposition}

\noindent\textbf{Proof.}
In light of Lemma~\ref{L:factor},
$\sum_\rho$ can also be thought of as the sum
over the non-trivial zeros of $L(s,\chi^k)$ for $k=1,\dots,\ell$
(counting multiplicities).
Observe that since $\ell$ is odd, we have $\chi(-1)=1$.
Now set $a=1/2$ and combine Lemmas
\ref{L:rlemma1}, \ref{L:rlemma2}, \ref{L:q2lemma3}, \ref{L:rlemma4}.
\qed

\vspace{1ex}
\noindent\textbf{Proof of Theorem~\ref{T:GRHr}.}
Define $x:=2.5(\ell-1)^2(\log f)^2$.
Since $f\geq 10^8$ and $\ell\geq 3$, we have $x>3393$.
By way of contradiction,
suppose that $\chi(n)\neq\omega$ for all $n<x$ with $(n, q_1 q_2)=1$.
Under this assumption we apply
Lemma~\ref{L:sumbound} with $u=q_1 q_2$, which gives
$$
  \sum_{\footnotesize{\begin{array}{cc}n<x\\[-1ex]\chi(n)=\omega\end{array}}}
  \Lambda(n)(n/x)^{1/2}\log(x/n)
  \leq 2(\log x)^2
\,.
$$

Combining the above with Proposition~\ref{P:explicit} and dividing by $\sqrt{x}$ yields:
$$
  \frac{\sqrt{x}}{9/4}
  +
  \frac{4}{\sqrt{x}}
  \leq
  \sum_\rho
  +
  \frac{2\ell(\log x)^2}{\sqrt{x}}
  +
  \frac{\log x}{x}
  \left[
  \frac{5}{2}\sum_\rho+\ell
  \left|\frac{\zeta'(2)}{\zeta(2)}\right|+\frac{25}{6}
  \right]
  +
  \frac{1}{x}
  \left[
  \sum_\rho
  +
  \frac{40}{9}
  \right]
$$

We note that Lemma~\ref{L:rzerosbound} is applicable in our situation;
indeed, our assumptions on $K$ imply that it is totally-real and, using
the conductor-discriminant formula, we see that $\Delta=f^{\ell-1}$.
By Lemma~\ref{L:rzerosbound},
we have
$$
  \sum_\rho\leq
  \frac{1}{2}
  \left[
  (\ell-1)\log f-
  2.23\,\ell
  +5.34
  \right]
  \,,
$$
and in particular,
$$
  \frac{1}{\sqrt{x}}\sum_\rho\leq\frac{1}{2\sqrt{2.5}}
  \,.
$$
We see
\begin{eqnarray*}
  \frac{\log x}{x}
  \left[
  \frac{25}{6}
  \right]
  +
  \frac{1}{x}
  \left[
  \sum_\rho
  +
  \frac{40}{9}
  \right]
   \;\leq\;
  \frac{1}{\sqrt{x}}
  \left(
        \frac{\log x}{\sqrt{x}}
  \cdot
  \frac{25}{6}
  +
  \frac{1}{2\sqrt{2.5}}
  +
  \frac{1}{\sqrt{x}}
  \frac{40}{9}
  \right)
  \;<\;
  \frac{4}{\sqrt{x}}
  \,,
\end{eqnarray*}
and therefore
$$
  \frac{\sqrt{x}}{9/4}
  \leq
  \sum_\rho
  +
  \frac{2\ell(\log x)^2}{\sqrt{x}}
  +
  \frac{\log x}{x}
  \left[
  \frac{5}{2}\sum_\rho+\ell
  \left|\frac{\zeta'(2)}{\zeta(2)}\right|
  \right]
  \,.
$$

We have
$$
  \frac{1}{\sqrt{x}}\left[
  \frac{5}{2}\sum_\rho
  +
  \ell\left|\frac{\zeta'(2)}{\zeta(2)}\right|
  \right]
  <
 0.82
$$
and
$$
\frac{\log x}{\sqrt{x}}
\leq
0.14
\,,
$$
which leads to
$$
  \frac{\sqrt{x}}{9/4}
  \leq
    \frac{1}{2}
  \left[
  (\ell-1)\log f-
  2.23\,\ell
  +5.34
  \right]
  +
  \frac{2\ell(\log x)^2}{\sqrt{x}}
  +
  0.12
  \,.
$$
Now we observe
\begin{eqnarray*}
  \frac{2\ell(\log x)^2}{\sqrt{x}}
  &\leq&
  2.27\ell
  \,.
\end{eqnarray*}

All together, we have
\begin{eqnarray*}
  \frac{\sqrt{x}}{9/4}
  &\leq&
      \frac{1}{2}
  \left[
  (\ell-1)\log f-
  2.23\,\ell
  +5.34
  \right]
  +
  2.27\ell
  +
  0.12
  \\
  &\leq&
  \frac{1}{2}(\ell-1)(\log f)+1.16\ell+2.79
\end{eqnarray*}
This leads to:
\begin{eqnarray*}
  \sqrt{x}
  &\leq&
  \frac{9}{8}(\ell-1)(\log f)
  +
  2.61\ell
  +
  6.28
  \\
  &\leq&
  1.51(\ell-1)\log f
\end{eqnarray*}
Squaring both sides yields
$$
  x\leq 2.3(\ell-1)^2(\log f)^2
  \,,
$$
a contradiction.
\qed


\section{GRH Bounds for Norm-Euclidean Fields}\label{S:GRH.DB}

In this section we prove Theorem~\ref{T:GRH}.
First we deal separately with the situation where $q_1$ is small.

\begin{theorem}\label{T:q1_100}
 Let $K$ be a Galois number field of odd prime degree $\ell$ and conductor $f$.
 Assume the GRH for $\zeta_K(s)$.
  Let $q_1$ denote the smallest rational prime which is inert in $K$.
  If $q_1<100$ and
$$
   5825(\ell-1)^2(\log f)^4<f
  \,,
$$
then $K$ is not norm-Euclidean.  
\end{theorem}

\noindent\textbf{Proof.}
Set $A=5825$.
One checks that our hypothesis implies $f\geq 10^9$ and $f>\ell^2$.
By Lemma~\ref{L:disc} we know that $f$ is a prime with $f\equiv 1\pmod{\ell}$.
We adopt the notation from the statement of Theorem~\ref{T:conditions}.
Since $q_1<100$, we have 
$q_1\leq 97$,
and by Theorems~\ref{T:GRHq2} and~\ref{T:GRHr}, we have
\begin{eqnarray}
  \label{E:est1}
  q_2&<&2.5(\log f)^2
  \,,
  \\
  \label{E:est2}
  r&<&2.5(\ell-1)^2(\log f)^2
  \,;
\end{eqnarray}
hence we have
\begin{eqnarray*}
  2.1\,q_1q_2 r\log q_1
  &<&
  (2.1)(97)(2.5)(\log f)^2(2.5)(\ell-1)^2(\log f)^2(\log 97)
  \\
  &<&
  A(\ell-1)^2(\log f)^4
  \,.
\end{eqnarray*}

If $q_1\neq 2,3,7$, then
it follows from Theorem~\ref{T:conditions}
that the condition given in our hypothesis is sufficient.
If $q_1=7$, then we observe that
\begin{eqnarray*}
  3\,q_1q_2 r\log q_1
  &<&
  (3)(7)(2.5)(\log f)^2(2.5)(\ell-1)^2(\log f)^2(\log 7)
  \\
  &<&
  A(\ell-1)^2(\log f)^4
  \,.
\end{eqnarray*}

Now we deal with the special case where $q_1=2$, $q_2=3$.
Our hypothesis gives
$$
  (\ell-1)<\frac{1}{\sqrt{A}}\,\frac{f^{1/2}}{(\log f)^2}
  \,.
$$
In order to use
Proposition~\ref{P:special}, we estimate
\begin{eqnarray*}
  72(\ell-1)f^{1/2}\log 4f+35
  &<&
  \frac{72}{\sqrt{A}}\,\frac{\log 4f}{(\log f)^2}\,f+35
  \\[1ex]
  &<&
  0.1 f+35
  \\
  &<&
  f
  \,;
\end{eqnarray*}
thus the proposition applies.
When $q_1=3$, $q_2=5$, we use a similar estimate to conclude that
$$
  507(\ell-1)f^{1/2}\log 9f+448<.4f+448<f
  \,,
$$
and hence Proposition~\ref{P:special} applies again.

The remaining cases fall under conditions (4) and (5) of Theorem~\ref{T:conditions}.
We will prove the bound
$$
  5\, q_2 r<f
  \,,
$$
which will deal with all remaining cases.
From the estimates (\ref{E:est1}) and (\ref{E:est2}) we have
\begin{eqnarray*}
  5\, q_2 r
  &<&
  32(\ell-1)^2(\log f)^4
  \\
  &<&
  A(\ell-1)^2(\log f)^4
  \\
  &<&
  f
  \,.
 \end{eqnarray*}
 This completes the proof.
 \qed

\vspace{1ex}

Applying the previous theorem with $\ell=3$ yields:
\begin{corollary}\label{C:myspecial}
 Let $K$ be a Galois cubic number field with conductor $f\geq 6\cdot 10^{9}$.
 Assume the GRH for $\zeta_K(s)$.
   Let $q_1$ denote the smallest rational prime which is inert in $K$.
  If $q_1<100$, then $K$ is not norm-Euclidean.
\end{corollary}

%

\noindent\textbf{Proof of Theorem~\ref{T:GRH}.}
One checks that our hypothesis implies
$f\geq 10^{10}$ and $f>\ell^2$.
We will apply Theorem~\ref{T:conditions} as in the proof
of Theorem~\ref{T:q1_100}.
Applying Theorems~\ref{T:bach}, \ref{T:GRHq2}, and~\ref{T:GRHr}, we have:
\begin{eqnarray}
  \nonumber
  q_1&<&(1.17\log f-6.3)^2\\
  \label{E:newq2}
  q_2&\leq&2.5(\log f)^2\\
  \label{E:newr}
  r&\leq&2.5(\ell-1)^2(\log f)^2
\end{eqnarray}
For the moment, we assume $q_1\neq 2,3,7$.
Combining everything, this gives
$$
  2.1\,q_1 q_2 r\log q_1<
  26.25(\ell-1)^2(1.17\log f-6.3)^2\log(1.17\log f-6.3)(\log f)^4
  \,.
$$
Hence a sufficient condition is:
\begin{equation}\label{E:ugly.condition}
  26.25(\ell-1)^2(1.17\log f-6.3)^2\log(1.17\log f-6.3)(\log f)^4\leq f
\end{equation}
Note that the condition given in our hypothesis implies (\ref{E:ugly.condition}).
To deal with the remaining cases of $q_1=2,3,7$, we note that 
(\ref{E:ugly.condition}) implies the condition given in the statement of
Theorem~\ref{T:q1_100}; hence (\ref{E:ugly.condition})
is sufficient in all cases.
\qed

%


\section{Galois Cubic Fields}\label{S:cubic}

Finally, we give the proof
of Theorem~\ref{T:cubic.GRH}.
Let $K$ be a norm-Euclidean Galois cubic field with conductor $f$ and discriminant $\Delta$
which is not
any of the $13$ fields listed in the statement of Theorem~\ref{T:cubic.GRH}.
In light of Theorem~\ref{T:cubic}, we may assume $f\geq 10^{10}$.
Moreover, Lemma~\ref{L:disc} allows us to conclude that $\Delta=f^2$ and that
$f$ is a prime with $f\equiv 1\pmod{3}$.
Using the
slightly complicated condition (\ref{E:ugly.condition})
in the proof of Theorem~\ref{T:GRH} and setting $\ell=3$ we
find that $f<7\cdot 10^{10}$.

It remains to deal with the cases where
$f$ lies in $(10^{10},\, 7\cdot 10^{10})$.
Let $\chi$ be a primitive cubic character modulo $f$,
and let $q_1$ denote the smallest prime such that $\chi(q_1)\neq 1$.
By Corollary~\ref{C:myspecial}, to show that $K$ is not norm-Euclidean,
assuming $f\in(10^{10},\, 7\cdot 10^{10})$,
it suffices to show $q_1<100$.  Using the method of character evaluation described 
in \S5.2 of~\cite{mcgown:euclidean}, we obtain the following lemma which completes
the proof of Theorem~\ref{T:cubic.GRH}.
\begin{lemma}
Suppose $f$ is a prime with $f\equiv 1\pmod{3}$.
Let $\chi$ be a cubic character modulo $f$,
and denote by $q_1$ the smallest prime with $\chi(q_1)\neq 1$.
If $f\leq 7\cdot 10^{10}$, then $q_1\leq 61$. 
\end{lemma}
The computation given in the above Lemma was carried out on an iMac
with a 3.06 GHz Intel Core 2 Duo processor and 4 GB of RAM, running Mac OS 10.6.
It took 8.4 days of CPU time to complete.


As an additional curiosity we have kept a list of record values of $q_1$.
That is, each time we encounter a value of $q_1$ which is strictly greater than all
previous values, we have outputted the values of $f$ and $q_1$.  Here are the results:

\vspace{1ex}
{\footnotesize
\begin{verbatim}
     Record: f=7, q1=2
     Record: f=31, q1=3
     Record: f=307, q1=5
     Record: f=643, q1=7
     Record: f=5113, q1=11
     Record: f=21787, q1=13
     Record: f=39199, q1=17
     Record: f=360007, q1=23
     Record: f=4775569, q1=29
     Record: f=10318249, q1=37
     Record: f=65139031, q1=41
     Record: f=387453811, q1=43
     Record: f=913900417, q1=47
     Record: f=2278522747, q1=53
     Record: f=2741702809, q1=59
     Record: f=25147657981, q1=61
\end{verbatim}
}
\section*{Acknowledgement}
This paper is based upon a portion of the author's Ph.D. dissertation under
the supervision of Professor Harold Stark.  The author would like to thank
Professor Stark for his invaluable guidance and support at all stages of this work.

\bibliographystyle{amsplain}


\bibliography{myrefs}  

\providecommand{\bysame}{\leavevmode\hbox to3em{\hrulefill}\thinspace}
\providecommand{\MR}{\relax\ifhmode\unskip\space\fi MR }
\providecommand{\MRhref}[2]{%
  \href{http://www.ams.org/mathscinet-getitem?mr=#1}{#2}
}
\providecommand{\href}[2]{#2}
\begin{thebibliography}{10}

\bibitem{ankeny}
N.~C. Ankeny, \emph{The least quadratic non residue}, Ann. of Math. (2)
  \textbf{55} (1952), 65--72.

\bibitem{bach:1990}
Eric Bach, \emph{Explicit bounds for primality testing and related problems},
  Math. Comp. \textbf{55} (1990), no.~191, 355--380.

\bibitem{davenport:mnt}
Harold Davenport, \emph{Multiplicative number theory}, third ed., Graduate
  Texts in Mathematics, vol.~74, Springer-Verlag, New York, 2000, Revised and
  with a preface by Hugh L. Montgomery.

\bibitem{godwin.smith:1993}
H.~J. Godwin and J.~R. Smith, \emph{On the {E}uclidean nature of four cyclic
  cubic fields}, Math. Comp. \textbf{60} (1993), no.~201, 421--423.

\bibitem{heilbronn:cubic}
H.~Heilbronn, \emph{On {E}uclid's algorithm in cubic self-conjugate fields},
  Proc. Cambridge Philos. Soc. \textbf{46} (1950), 377--382.

\bibitem{murty:2009}
Yasutaka Ihara, V.~Kumar Murty, and Mahoro Shimura, \emph{On the logarithmic
  derivatives of {D}irichlet {$L$}-functions at {$s=1$}}, Acta Arith.
  \textbf{137} (2009), no.~3, 253--276.

\bibitem{lagarias:1977}
J.~C. Lagarias and A.~M. Odlyzko, \emph{Effective versions of the {C}hebotarev
  density theorem}, Algebraic number fields: {$L$}-functions and {G}alois
  properties ({P}roc. {S}ympos., {U}niv. {D}urham, {D}urham, 1975), Academic
  Press, London, 1977, pp.~409--464.

\bibitem{landau:1}
Edmund Landau, \emph{Zur {T}heorie der {H}eckeschen {Z}etafunktionen, welche
  komplexen {C}harakteren entsprechen}, Math. Z. \textbf{4} (1919), no.~1-2,
  152--162.

\bibitem{lemmermeyer:euclidean}
Franz Lemmermeyer, \emph{The {E}uclidean algorithm in algebraic number fields},
  Exposition. Math. \textbf{13} (1995), no.~5, 385--416.

\bibitem{mcgown:euclidean}
Kevin~J. McGown, \emph{Norm-{E}uclidean {G}alois fields},  (to appear).

\bibitem{montgomery:book}
Hugh~L. Montgomery, \emph{Topics in multiplicative number theory}, Lecture
  Notes in Mathematics, Vol. 227, Springer-Verlag, Berlin, 1971.

\bibitem{narkiewicz:book}
W{\l}adys{\l}aw Narkiewicz, \emph{Elementary and analytic theory of algebraic
  numbers}, second ed., Springer-Verlag, Berlin, 1990.

\bibitem{odlyzko:1990}
A.~M. Odlyzko, \emph{Bounds for discriminants and related estimates for class
  numbers, regulators and zeros of zeta functions: a survey of recent results},
  S\'em. Th\'eor. Nombres Bordeaux (2) \textbf{2} (1990), no.~1, 119--141.

\bibitem{poitou:1}
Georges Poitou, \emph{Minorations de discriminants (d'apr\`es {A}. {M}.
  {O}dlyzko)}, S\'eminaire {B}ourbaki, {V}ol. 1975/76 28\`eme ann\'ee, {E}xp.
  {N}o. 479, Springer, Berlin, 1977, pp.~136--153. Lecture Notes in Math., Vol.
  567.

\bibitem{poitou:2}
\bysame, \emph{Sur les petits discriminants}, S\'eminaire
  {D}elange-{P}isot-{P}oitou, 18e ann\'ee: (1976/77), {T}h\'eorie des nombres,
  {F}asc. 1 ({F}rench), Secr\'etariat Math., Paris, 1977, pp.~Exp. No. 6, 18.

\bibitem{stark:1974}
H.~M. Stark, \emph{Some effective cases of the {B}rauer-{S}iegel theorem},
  Invent. Math. \textbf{23} (1974), 135--152.

\bibitem{stark:1975}
\bysame, \emph{The analytic theory of algebraic numbers}, Bull. Amer. Math.
  Soc. \textbf{81} (1975), no.~6, 961--972.

\bibitem{weil:ef.1}
Andr{\'e} Weil, \emph{Sur les ``formules explicites'' de la th\'eorie des
  nombres premiers}, Comm. S\'em. Math. Univ. Lund [Medd. Lunds Univ. Mat.
  Sem.] \textbf{1952} (1952), no.~Tome Supplementaire, 252--265.

\bibitem{weil:ef.2}
\bysame, \emph{Sur les formules explicites de la th\'eorie des nombres}, Izv.
  Akad. Nauk SSSR Ser. Mat. \textbf{36} (1972), 3--18.

\end{thebibliography}

\end{document}